\documentclass[12pt]{article}
\usepackage[dvipsnames,usenames]{color}
\usepackage{amsmath}
\usepackage{amsfonts}
\usepackage{amssymb}
\def\h{\hbox}

\newcommand{\CC}{{\mathbb C}}
\newcommand{\RR}{{\mathbb R}}
\newcommand{\HH}{{\mathbb H}}
\newcommand{\BB}{{\mathbb B}}
\newcommand{\PP}{{\mathbb P}}
\newcommand{\EE}{{\mathbb E}}

\newcommand{\ra}{\rightarrow}
\newcommand{\ov}{\overline}
\newcommand{\p}{\partial}
\newcommand{\w}{\widetilde}

                     \newtheorem{theorem}{Theorem}[section]
\newtheorem{lemma}[theorem]{Lemma}

\begin{document}
\author{Shanyu Ji and Yuan Yuan}
\title{\bf Flatness of CR Submanifolds in a Sphere}
\maketitle


\bigskip
\centerline{\it Dedicated to Professor Yang, Lo in the Occasion of his 70th Birthday}
\bigskip

\section{Introduction}
The Cartan-Janet theorem asserted that for any analytic Riemannian
manifold $(M^n, g)$, there exist local isometric embeddings of $M^n$
into Euclidean space $\EE^N$ as $N$ is sufficiently large.
The CR analogue of Cartan-Janet theorem is not true
in general. In fact, Forstneric \cite{Fo86} and Faran \cite{Fa88}
proved the existence of real analytic strictly pseudoconvex
hypersurfaces $M^{2n+1}\subset \CC^{n+1}$ which do not admit any
germ of holomorphic mapping taking $M$ into sphere $\p\BB^{N+1}$ for
any $N$.

There are recent progress on CR submanifolds in sphere $\p\BB^{N+1}$.
Zaitsev \cite{Za08} constructed explicit examples for the Forstneric and Faran phenomenon
above. Ebenfelt, Huang and Zaitsev \cite{EHZ04}
proved rigidity of CR embeddings of general $M^{2n+1}$ into spheres with
CR co-dimension $<\frac{n}{2}$, which generalizes a result of Webster that was for the case
of co-dimension 1 \cite{We79}.
S.-Y. Kim and J.-W. Oh \cite{KO06} gave a necessary and sufficient condition for local
embeddability into a sphere $\p\BB^{N+1}$
of a generic strictly pseudoconvex psuedohermitian CR manifold $(M^{2n+1}, \theta)$ in terms of its Chern-Moser curvature tensors and their derivatives.

In Euclidean geometry, for a real submanifold $M^n\subset
\EE^{n+a}$, $M$ is a piece of $\EE^n$  if and only if
its second fundamental form $II_M\equiv 0$.
In projective geometry, for a complex submanifold $M^n\subset
\CC\PP^{n+a}$,
$M$ is a piece of $\CC\PP^n$ if and only if
its projective second fundamental form $II_M\equiv 0$ (c.f. \cite{IL03}, p.81).
In CR geometry,  we prove the CR analogue of this fact in this paper as follows:

\begin{theorem}
Let $H: M' \ra \p\BB^{N+1}$ be a smooth CR-embedding of a strictly pseudoconvex CR real hypersurface
$M' \subset \CC^{n+1}$. Denote $M:=H(M')$.
If its CR second fundamental form $II_M\equiv 0$, then $M\subset F(\p\BB^{n+1})\subset \p\BB^{N+1}$ where
$F:\BB^{n+1}\ra \BB^{N+1}$ is a certain linear fractional proper holomorphic map.
\end{theorem}

Previously, it was proved by P. Ebenfelt, X. Huang and D. Zaitsev (\cite{EHZ04},
corollary 5.5), under the above same hypothese, that $M'$ and hence $M$ are
locally CR-equivalent to the unit sphere
$\p\BB^{n+1}$ in $\CC^{n+1}$.

There are several definitions of the CR second fundamental forms $II_M$ of $M$ (see
Section 3, 4, 5, and 6). The result in \cite{EHZ04} used Definition 1 or 2. However, to prove Theorem
1.1, we need to use Definitions 3 and 4. We'll prove in Section 4 that $II_M\equiv 0$
by any one of the four definitions will imply $II_M\equiv 0$ for all other three definitions.
One of the
ingredients for our proof of Theorem 1.1 is the result of Ebenfelt-Huang-Zaitsev \cite{EHZ04} so that
$M$ can be regarded as the image of a rational CR map $F:\p\HH^{n+1}\ra M \subset \p\HH^{N+1}$.
Another ingredient is a theorem of Huang (\cite{Hu99})
that such a map $F$ is linear if and only if its geometric rank $\kappa_0$ is zero.
The third one is a result from \cite{HJY09} about a special lift for maps between spheres.

\bigskip

{\bf Acknowledgments} We would like to thank
Professor Xiaojun Huang for the constant encouragement and support.
The second author is also grateful to Wanke Yin and Yuan Zhang for helpful
discussions.

\section{Preliminaries}


\noindent $\bullet$ {\bf Maps between balls }
We denote by $Prop(\BB^n, \BB^N)$ the space of all proper holomorphic maps from the unit ball $\BB^n\subset \CC^n$
to $\BB^N$, denote by $Prop_k(\BB^n, \BB^N)$ the space $Prop(\BB^n, \BB^N)\cap C^k(\ov{\BB^n})$, and
denote by $Rat(\BB^n, \BB^N)$ the space $Prop(\BB^n, \BB^N)\cap \{rational\ maps\}$.
We say that $F$ and $G\in Prop(\BB^{n}, \BB^{N})$ are
{\it equivalent} if there are automorphisms $\sigma\in Aut(\BB^{n})$
and $\tau\in Aut(\BB^{N})$ such that $F=\tau\circ G\circ\sigma$.

Write $\HH^n:=\{(z,w)\in \CC^{n-1}\times \CC: \
\hbox{Im}(w)>|z|^2\}$ for the Siegel upper-half space. Similarly, we
can define the space $Prop(\HH^n, \HH^N)$, $Prop_k(\HH^n, \HH^N)$
and $Rat(\HH^n, \HH^N)$ similarly. By the Cayley transformation
$\rho_n: \HH^n\to \BB^n$, $\rho_n(z,w)=(\frac{2z}{1-iw},\
\frac{1+iw}{1-iw})$, we can identify a map $F\in Prop_k(\BB^n,
\BB^N)$ or $Rat(\BB^n, \BB^N)$ with $\rho^{-1}_N\circ F\circ\rho_n$
in the space $Prop_k(\HH^n, \HH^N)$ or $Rat(\HH^n, \HH^N)$,
respectively. We say that $F$ and $G\in Prop(\HH^{n}, \HH^{N})$ are
{\it equivalent} if there are automorphisms $\sigma\in Aut(\HH^{n})$
and $\tau\in Aut(\HH^{N})$ such that $F=\tau\circ G\circ\sigma$.


We denote by $\p\HH^n=
\{(z,w)\in \CC^{n-1}\times \CC:\ \hbox{Im}(w)=|z|^2\}$ for the Heisenberg hypersurface.
For any map $F\in Prop_2(\HH^n, \HH^N)$, by restricting on $\p\HH^n$, we can regard $F$ as a
$C^2$ CR map from $\p\HH^n$ to $\p\HH^N$, and we denote it as $F\in Prop_2(\p\HH^n, \p\HH^N)$.
We say that $F$ and $G\in Prop_2(\p\HH^{n}, \p\HH^{N})$ are
{\it equivalent} if there are automorphisms $\sigma\in Aut(\p\HH^{n})=Aut(\HH^n)$
and $\tau\in Aut(\p\HH^{N})=Aut(\HH^N)$ such that $F=\tau\circ G\circ\sigma$.

We can parametrize $\partial \mathbb{H}^n$ by $(z,\overline{z},u)$ through the
map $(z,\overline{z},u)\to (z,u+i|z|^2)$. In what follows, we will assign
the weight of $z$ and $u$ to be   $1$ and $2$, respectively.
For a non-negative integer $m$, a function
$h(z,\overline{z},u)$ defined over a small ball  $U$
of $0$ in $\partial \mathbb{H}^n$  is said to be of quantity $o_{wt}(m)$
if
$\frac{h(tz,t\overline{z},t^2u)}{|t|^{m}}\to 0$ uniformly for $(z,u)$ on
any
             compact subset of $U$ as $t(\in \mathbb{R})\to 0$.

\bigskip

\noindent $\bullet$ {\bf Partial normalization  of $F$} Let
$F=(f,\phi,g)=(\widetilde{f}, g)= (f_1,\cdots,f_{n-1},
\phi_1,\cdots, \phi_{N-n},g)$ be a non-constant map in
$Prop_2(\p\HH^n, \p\HH^N)$ with $F(0)=0$. For each $p\in
\p\mathbb{H}^n$, we write $\sigma^0_p \in \hbox{Aut}(\mathbb{H}^n)$
with $\sigma^0_p(0)=p$ and
$\tau^F_p\in\hbox{Aut}(\mathbb{H}^N)$ with
$\tau^F_p(F(p))=0$ for the maps

\begin{eqnarray}
\label{sigma}
&&\sigma^0_p(z,w)=(z+z_0, w+w_0+2i \langle z,\overline{z_0}
\rangle),\\
&&\tau^F_p(z^*,w^*)=(z^*-\widetilde{f}(z_0,w_0),w^*-\overline{g(z_0,w_0)}-
2i \langle z^*,\overline{\widetilde{f}(z_0,w_0)} \rangle ). \label{tau}
\end{eqnarray}

$F$ is equivalent
to $F_p=\tau^F_p\circ F\circ \sigma^0_p=(f_p,\phi_p,g_p). $ Notice that
$F_0=F$ and $F_p(0)=0$.
The following is basic for the understanding of the geometric properties
of
$F$.

\bigskip
\begin{lemma}
\label{**normalform}
([$\S 2$, Lemma 5.3, Hu99], [Lemma 2.0, Hu03]):
Let $F$ be a non-constant map in $Prop_2(\p\HH^n, \p\HH^N)$, $2\le n\le
N$ with $F(0)=0$.
For each $p\in \partial \mathbb{H}^n$, there is an automorphism
$\tau^{**}_p\in Aut_0(\HH^N)$ such that
$F_{p}^{**}:=\tau^{**}_p\circ F_p$ satisfies
the following  normalization:
\begin{equation}
f^{**}_{p}=z+{\frac{i}{ 2}}a^{**(1)}_{p}(z)w+o_{wt}(3),\ \phi_p^{**}
={\phi_p^{**}}^{(2)}(z)+o_{wt}(2),
\ g^{**}_{p}=w+o_{wt}(4),\
\label{eqn: **formal}
\end{equation}
$$\langle \overline{z}, a_{p}^{**(1)}(z)\rangle
|z|^2=|{\phi_p^{**}}^{(2)}(z)|^2.$$
\end{lemma}

\bigskip
Let $\mathcal{A}(p)=-2i(\frac{\partial^2(f_p)^{\ast\ast}
             _l}{\partial z_j\partial w}|_0)_{1\leq j,l\leq n-1}$. We
call the rank of $\mathcal{A}(p)$, which we denote by
$Rk_F(p)$, the {\it geometric rank} of $F$ at $p$.
$Rk_F(p)$ depends only on $p$ and $F$, and is a lower
semi-continuous function on $p$.
We define the {\it geometric rank} of $F$ to be
$\kappa_0(F)=max_{p\in \partial\mathbb{H}^n} Rk_F(p)$.
Notice that we always have $0\le \kappa_0\le n-1$.
We define the geometric rank of $F\in
\h{Prop}_2(\mathbb{B}^n, \mathbb{B}^N)$ to be the one for
the map $\rho_N^{-1}\circ F \circ \rho_n\in \h{Prop}_2(\mathbb{H}^n,
\mathbb{H}^N)$.


\begin{lemma}
(ct. [Hu99], theorem 4.3)
\label{rank 0}
$F\in Prop_2(\BB^n, \BB^N)$ has geometric rank $0$ if and only if $F$ is
equivalent to a linear map.
\end{lemma}

Denote by ${\cal S}_{0}=\{(j,l): 1\le j\leq \kappa_0, 1\leq l\leq (n-1),
j\leq l\}$ and write
${\cal S}:=\{(j,l):\ (j,l)\in {\cal S}_0,\ \h{or}\ \ j=\kappa_0+1,
l\in \{\kappa_0+1,\cdots, \kappa_0+N-n-\frac{(2n-\kappa_0-1)\kappa_0}{
2}\}\}$.

\bigskip

\begin{lemma}
\label{***normalform}
([Lemma 3.2, Hu03]):
Let $F$ be a $C^2$-smooth CR map from an open
piece $M\subset \partial\mathbb{H}^n$ into
$\partial\mathbb{H}^{N}$ with $F(0)=0$ and $Rk_F(0)=\kappa_0$.
Let $P(n,\kappa_0)=\frac{\kappa_0(2n-\kappa_0-1)}{2}$. Then $N\ge n+P
(n,\kappa_0)$
and  there are $\sigma\in \h{Aut}_0(\partial\mathbb{H}^n)$ and $\tau
\in\h{Aut}_0
(\partial \mathbb{H}^N)$
such that $F^{***}_p=\tau\circ F\circ \sigma:=(f,\phi,g)$ satisfies the
following normalization conditions:
\begin{equation}
\left\{
\begin{aligned}
f_j=&z_j+\frac{i\mu_j}{2}z_jw+o_{wt}(3),\ \ \frac{\p^2 f_j}{\p w^2}(0)=0,\
j=1
\cdots,\kappa_0,\ \mu_j>0,\cr
f_j=&z_j+o_{wt}(3), \ \ j=\kappa_0+1,\cdots, n-1\cr
             g=&w+o_{wt}(4),\cr
\phi_{jl}=&\mu_{jl}z_jz_l+o_{wt}(2), \ \h{where}\ \ (j,l)\in {\cal S}\
                      \h{with}\ \ \mu_{jl}>0 \ \h{for}\ \  (j,l)\in {\cal
S}_0\cr
                         &\h{and}\ \mu_{jl}=0\ \h{otherwise}
\end{aligned}
\right.
\end{equation}
where $\mu_{jl}=\sqrt{\mu_j+\mu_l}$ for $j,l\le
\kappa_0$ $j\not=l$,
$\mu_{jl}=\sqrt {\mu_j}$ if $j\le \kappa_0$ and $l>\kappa_0$ or if
$j=l\le \kappa_0$.
\end{lemma}

\bigskip

\noindent$\bullet$ {\bf Pseudohermitian metric and Webster connection}\ \ \ \
Let $M$ be a  $C^2$ smooth real hypersurface in $\CC^{n+1}$.
We denote by $T^c M=TM\cap i TM\subset TM$ its {\it maximal complex tangent bundle} with
the complex structure $J: T^c M \ra T^c M$. Here $J(\frac{\p}{\p x_j})=\frac{\p}{\p y_j}$ and
$J(\frac{\p}{\p y_j})=-\frac{\p}{\p x_j}$ in terms of holomorphic coordinates.
We denote by
${\cal V}=T^{0,1}M=\{X+i J X\ |\ X\in T^c M\}\subset \CC TM:=TM\otimes \CC$ the
{\it CR bundle}.
We also denote $T^{1,0}M=\ov{\cal V}$. All $T^c M$, ${\cal V}$ and $\ov{\cal V}$ are
complex rank $n$ vector bundles.

Write $T^0 M:=(T^{1,0}M\oplus T^{0,1}M)^\perp \subset \CC T^*M$ for its
rank one subbundle.  Write $T'M:={T^{0,1}}^\perp \subset \CC T^* M$
for its rank $n+1$ {\it holomorphic or (1,0) cotangent bundle} of
$M$. Here $T^0\subset T'M$.

A real nonvanishing $1$-form  $\theta$ over $M$ is called  a
{\it contact form} if $\theta\wedge (d \theta)^n\not=0$.
Let $M$ be as above given by a
defining function $r$. Then the $1$-form $\theta=i \p r$ is a contact form of $M$.

We say that $(M, \theta)$ is {\it strictly pseudoconvex} if the Levi-form $L_\theta$
is positive definite for all $z\in M$. Here the {\it Levi-form} $L_\theta$ with respect to
$\theta$ is defined by
\[
L_\theta(\vec u, \ov{\vec v}):=-id\theta(\vec u \wedge \ov{\vec v})
,\ \ \ \forall \vec u,
\vec v\in T^{1,0}_p(M),\ \forall p\in M.
\]

Associated with a contact form $\theta$ one has the Reeb
vector field $R_\theta$, defined by the equations: (i) $d \theta(R_\theta, \cdot) \equiv 0$,\
(ii) $\theta (R_\theta) \equiv 1$.\ As a skew-symmetric form of maximal rank $2n$,
the form $d\theta|_{T_p M}$ has a 1-
dimensional kernel for each $p \in  M^{2n+1}$. Hence equation (i) defines a unique
line field $\langle R_\theta \rangle $ on $M$. The contact condition $\theta \wedge (d\theta)^n \not= 0 $
implies that $\theta $ is non-trivial on that line field, so the unique real vector
field is defined by the normalization condition (ii).

According Tanaka \cite{T75} and Webester \cite{We78},
$(M, \theta)$ is called a {\it strictly pseudoconvex pseudohermitian manifold}
if there are $n$ complex $1$-forms $\theta^\alpha$ so that $\{\theta^1, ..., \theta^n\}$
forms a local basis for holomorphic cotangent bundle $H^*(M)$ and
\begin{equation}
\label{CR metric}
d\theta=i\sum^n_{\alpha, \beta=1} h_{\alpha \ov{\beta}} \theta^\alpha \wedge
\theta^{\ov{\beta}}
\end{equation}
where $(h_{\alpha \ov{\beta}})$, called the {\it Levi form matrix}, is positive definite.
Such $\theta^\alpha$ may not be unique. Following Webster (1978), a coframe
$(\theta, \theta^\alpha)$ is called
{\it admissible} if (\ref{CR metric}) holds. The admissible coframes are determined up to
transformations
$\w \theta^\alpha=u^\alpha_\beta \theta^\beta$ where $(u^\alpha_\beta)\in GL(\CC^n)$.


\begin{theorem}
(Webster, 1978)
Let $(M^{2n+1}, \theta)$ be a strictly pseudoconvex pseudohermitian manifold and
let $\theta^j$ be as in (\ref{CR metric}). Then there are unique way to write
\begin{equation}
\label{CR structure eq}
d\theta^\alpha = \sum^n_{\gamma=1} \theta^\gamma \wedge \omega^\alpha_\gamma + \theta
\wedge \tau^\alpha,
\end{equation}
where $\tau^\alpha$ are $(0,1)$-forms over $M$ that are linear combination of
$\theta^{\ov\alpha}
=\ov{\theta^\alpha}$,
and $\omega^\beta_\alpha$ are $1$-forms over $M$ such that
\begin{equation}
\label{CR condition}
0=d h_{\alpha \ov\beta}-h_{\gamma \ov\beta} \omega^\gamma_\alpha-h_{\alpha\ov\gamma}
\omega^{\ov\gamma}_{\ov\beta}.
\end{equation}
\end{theorem}

\bigskip

We may denote $\omega_{\alpha\ov\beta}=h_{\gamma\ov\beta}\omega^\gamma_\alpha$ and
$\ov{\omega_{\beta\ov\alpha}}=h_{\alpha\ov\gamma}\omega^{\ov\gamma}_{\ov\beta}$.
In particular, if
\begin{equation}
\label{orthonormal}
h_{\alpha\beta}=\delta_{\alpha \beta},
\end{equation}
the identity in (\ref{CR condition})
becomes $0=-\omega_{\alpha\ov\beta}-\ov{\omega_{\beta\ov\alpha}}$,
i.e.,
\begin{equation}
\label{CR condition 2}
0=\omega^\beta_\alpha+\omega^{\ov\alpha}_{\ov \beta}.
\end{equation}
The condition on $\tau^\beta$ means:
\begin{equation}
\tau^\beta=A^\beta_{\ \ov\nu} \theta^{\ov{\nu}}, \ \ A^{\alpha\beta}=A^{\beta\alpha},
\end{equation}
which holds automatically. The curvature is given by
\begin{equation}
d \omega^{\ \beta}_\alpha-\omega^{\ \gamma}_\alpha \wedge \omega^{\ \beta}_\gamma
 = R^{\ \beta}_{\alpha\ \ \mu\ov\nu} \theta^\mu\wedge \theta^{\ov \nu}
 + W^{\ \beta}_{\alpha\ \ \mu}\theta^\mu\wedge \theta-W^\beta_{\ \ \alpha \ov \nu}
\theta^{\ov \nu}
\wedge \theta+i\theta_\alpha\wedge \tau^\beta-i\tau_\alpha\wedge \theta^\beta
\end{equation}
where the functions $R^{\ \beta}_{\alpha\ \ \mu\ov\nu}$ and $W^{\ \beta}_{\alpha\ \ \mu}$
represent the {\it pseudohermitian curvature} of $(M, \theta)$.

\bigskip

\section{CR second fundamental forms ----- Definition 1}
We are going to survey four definitions of the CR
second fundamental forms $II_M$ of $M$ in $\p\HH^{N+1}$. We start with Definition 1 which is the intrinsic
one in terms of a coframe.



\begin{lemma}
\label{lemma EHZ04} (\cite{EHZ04}, corollary 4.2) Let $M$ and $\w M$
be strictly pseudoconvex CR-manifolds of dimensions $2n+1$ and $2\w
n+1$ respectively, and of CR dimensions $n$ and $\w n$ respectively.
Let $F: M\ra \w M$ be a smooth CR-embedding. If $(\theta,
\theta^\alpha)$ is a admissible coframe on $M$, then  in a
neighborhood of a point $\w p\in F(M)$ in $\w M$ there exists an
admissible coframe $(\w \theta, \w \theta^A) = (\w
\theta, \w \theta^\alpha, \w \theta^\mu)$ on $\w M$ with $F^*(\w
\theta, \w \theta^\alpha, \w \theta^\mu)=(\theta, \theta^\alpha,
0)$. In particular, the Reeb vector field $\w R$ is tangent to
$F(M)$. If we choose the Levi form matrix of $M$ such that the
functions $h_{\alpha\ov{\beta}}$ in (\ref{CR metric}) with respect
to $(\theta, \theta^\alpha)$ to be $(\delta_{\alpha \ov\beta})$,
then $(\w \theta, \w \theta^A)$ can be chosen such that the Levi
form matrix of $\w M$ relative to it is also $(\delta_{A \ov B})$.
With this additional property, the coframe $(\w \theta, \w
\theta^A)$ is uniquely determined along $M$ up to unitary
transformations in $U(n)\times U(\w n-n)$.
\end{lemma}

If $(\theta, \theta^\alpha)$ and $(\w \theta, \w \theta^A)$ are as
above such that the condition on the Levi form matrices in Lemma
\ref{lemma EHZ04} are satisfied, we say that the coframe $(\w
\theta, \w \theta^A)$ is {\it adapted} to the coframe $(\theta,
\theta^\alpha)$. In this case, by (\ref{CR condition 2}), we have
$\theta=F^*\w\theta$, $\theta^\alpha=F^* \w\theta^\alpha$, and
\begin{equation*}
d\theta^\alpha = \sum^n_{\gamma=1} \theta^\gamma \wedge \omega^\alpha_\gamma + \theta
\wedge \tau^\alpha,\ \
0=\ \omega_\alpha^\beta+\omega_{\ov\beta}^{\ov\alpha},\ \ \ \forall 1\le \alpha, \beta\le n,
\end{equation*}
and
\begin{equation*}
d\w \theta^A = \sum^{\w n}_{B=1} \w\theta^C \wedge \w\omega^A_C + \w\theta
\wedge \w\tau^A,\ \
0=\ \w\omega_A^B+\w\omega_{\ov B}^{\ov A},\ \ \ \ \ \forall 1\le A, B\le N.
\end{equation*}
For simplicity, we may denote $F^* \w\omega^A_B$ by $\omega^A_B$.
We also denote $F^*\w\omega_{A \ov B}$ by $\omega_{A \ov B}$ where $\omega_{A \ov B}=
\omega_A^B$.

Write
$\omega^{\ \mu}_\alpha=\omega^{\ \mu}_{\alpha\ \beta}\theta^\beta$.
The matrix of $(\omega_{\alpha\ \ \beta}^{\ \mu})$, $1\le \alpha, \beta\le n$, $n+1\le \mu\le \hat n$,
defines the
{\it CR second fundamental form} of $M$. It was used in \cite{We79} and \cite{Fa90}.

\bigskip

\section{CR second fundamental forms ----- Definition 2}

Definition 2 introduced in \cite{EHZ04} will be the
extrinsic one in terms of defining function.

Let
$F: M \ra \w M$ be a smooth CR-embedding between $M\subset
\CC^{n+1}$ and $\w M\subset \CC^{N+1}$ where $M$ and $\w M$ are real
strictly pseudoconvex hypersurfaces of dimensions $2n+1$ and $2\w
n+1$, and CR dimensions $n$ and $\w n$, respectively. Let $p\in M$
and $\w p=F(p)\in \w M$ be points. Let $\w\rho$ be a local defining
function for $\w M$ near the point $\w p$. Let
\[
E_k(p):=span_\CC\{L^{\bar{J}} (\w\rho_{Z'}\circ F)(p)\ |\ J\in
(Z_+)^n, 0\le |J| \le k \}\subset T^{1,0}_{\w p} \CC^{N+1},
\]
where $\w\rho_{Z'}:=\p \w\rho$ is the complex gradient (i.e.,
represented by vectors in $\CC^{N+1}$ in some local coordinate
system $Z'$ near $\w p$). Here we use multi-index notation $L^{\ov
J}=L_1^{\ov{J_1}} \cdots L_n^{\ov{J_n}}$ and $|J|=J_1+...+J_n$.
It was shown in \cite{La01} that $E_k(p)$ is
independent of the choice of local defining function $\w \rho$,
coordinates $Z'$ and the choice of basis of the CR vector fields
$L_{\ov 1}, ..., L_{\ov n}$.

The {\it CR second fundamental form} $II_M$ of $M$ is defined by (cf. [EHZ04], $\S 2$)
\begin{equation}
\label{II M 2} II_M(X_p, Y_p):=\ov{\pi\big(XY(\w\rho_{\ov Z'}\circ
f)(p)\big)} \in \ov{T'_{\tilde p} \w M / E_1(p)}
\end{equation}
where $\w\rho_{\ov Z'}=\ov\p\w\rho$ is represented by vectors in
$\CC^{N+1}$ in some local coordinate system $Z'$ near $\w p$, $X, Y$
are any $(1,0)$ vector fields on $M$ extending given vectors $X_p,
Y_p\in T^{1,0}_p(M)$, and $\pi: T'_{\w p} \w M \ra T'_{\w p} \w
M/E_1(p)$ is the projection map.

Since $\w M$ and $M$ are strictly pseudoconvex, the Levi form of $\w M$ (at $\w p$) with respect to
$\w\rho$ defines an isomorphism
\[
\ov{T'_{\w p} \w M/E_1(p)}\cong T^{1,0}_{\w p}\w M/F_*(T^{1,0}_p M)
\]
and the CR second fundamental form can be viewed as an $\CC$-linear symmetric form
\begin{equation}
\label{bilinear form}
II_{M, p}: T^{1,0}_p M \times T^{1,0}_p M \ra T^{1,0}_{\w p} \w M / F_*(T^{1,0}_p M)
\end{equation}
that does not depend on the choice of $\w\rho$ (cf.[EHZ04], $\S 2$).

\bigskip

The relation between Definition 1 and Definition 2
was discussed in \cite{EHZ04}. Let $(M, \w M)$, $(\theta,
\theta^\alpha), (\w\theta, \w\theta^A)$ be as in Lemma 3.1,
and we abuse the structure bundle $(\theta,
\theta^ \alpha)$ on $M$ with the structure bundle $(\w \theta,
\w \theta ^\alpha)$ on $\w M$. We can
choose a defining function $\w\rho$ of $\w M$ near a point $\w
p=F(p) \in \w M$ where $p\in M$ such that $\theta=i \ov\p \w\rho$ on
$\w M$., i.e., in local coordinates $Z'$ in $\CC^{N+1}$, we have
\[
\theta=i \sum^{N+1}_{k=1} \frac{\p \w\rho}{\p \ov{Z'_k}} d \ov{Z_k'},
\]
where we pull back the forms $d \ov{Z'_1}, ..., d \ov{Z'_{N+1}}$ to
$\w M$. Then we consider the coframe $(\theta, \theta^\alpha)=(F^*
\w\theta, F^*\w\theta^\alpha)$ on $M$ near $p$ with $F(p)=\w p$. We
take its dual frame $(T, L_A)$ of $(\theta, \theta^A)$ and have

\begin{equation}
\label{2 CR SMF}
L_\beta(\w\rho_{\ov Z'} \circ F)
=-i L_\beta \lrcorner d \theta = g_{\beta \ov C}\theta^{\ov C} =
g_{\beta \ov\gamma} \theta^{\ov \gamma}.
\end{equation}

Here we used the definition of the construction, (\ref{CR metric})
and the dual relationship $\langle L_\beta,
\theta^\alpha\rangle=\delta_{\beta}^\alpha$ and also notice that
$g_{\beta \ov\gamma}=\delta_{\beta \gamma}$. Applying $L_\alpha$ to
both sides of (\ref{2 CR SMF}), we obtain
\[
L_\alpha L_\beta(\w\rho_{\ov Z'}\circ F)=g_{\beta \ov\gamma}
L_\alpha \lrcorner d\theta^{\ov \gamma} =\omega_{\alpha \ov \mu
\beta}\theta^{\ov \mu}\ \ \ mod(\theta, \theta^{\ov \alpha})
\]
which implies
\begin{equation}
\label{IIM b} II_M(L_\alpha, L_\beta)=\omega^{\ \mu}_{\alpha\ \beta}
L_\mu,\ \ \ n+1\le \mu \le N.
\end{equation}
This identity gives the equivalent relation of the intrinsic and
extrinsic definitions of $II_M$. Notice that we need a right choice
of $(\theta, \theta^\alpha)$, $(T, L_A)$ and $\w\rho$.


By using $(\omega^{\ b}_{\alpha\ \beta})$ and (\ref{IIM b}), as in (\ref{bilinear form}),
we can also define
\begin{equation}
\label{bilinear form 2}
II_{M, p}: T^{1,0}_p M \times T^{1,0}_p M \ra T^{1,0}_{\w p} \w M / F_*(T^{1,0}_p M)
\end{equation}
which is independent of the choice of the adapted coframe $(\theta,
\theta^A)$ in case $\w M$ is locally CR embeddable in $\CC^{N+1}$
(cf. [EHZ04], $\S$ 4).

\bigskip
\section{CR second fundamental forms ----- Definition 3}

\

Definition 3 will be the one as a tensor with respect to the group $GL^Q(\CC^{N+2})$.

\bigskip

\noindent {\bf The bundle  $GL^Q(\CC^{N+2})$ over $\p\HH^{N+1}$}\ \ \ \
We consider a real hypersurface $Q$ in $\CC^{N+2}$ defined by the homogeneous equation
\begin{equation}
\label{3 sphere 2}
\langle Z, Z\rangle : = \sum_A Z^A \ov{Z^A} + \frac{i}{2}(\ov{Z^0}Z^{N+1}
- Z^0
\ov{Z^{N+1}})=0,
\end{equation}
where $Z=(Z^0, Z^A, Z^{N+1})^t\in \CC^{N+2}$. Let
\begin{equation}
\label{pi0}
\pi_0: \CC^{N+2}-\{0\} \ra \CC\PP^{N+1},\ \  (z_0, ...., z_{N+1}) \mapsto [z_0: ... : z_{N+1}],
\end{equation}
be the standard projection. For any point $x\in \CC\PP^{N+1}$, $\pi_0^{-1}(x)$ is a complex line
in $\CC^{N+2}-\{0\}$. For any point $v\in \CC^{N+2}-\{0\}$, $\pi_0(v)\in \CC\PP^{N+1}$ is a point.
The image $\pi_0(Q-\{0\})$ is the Heisenberg hypersurface $\p\HH^{N+1}\subset \CC\PP^{N+1}$.

For any element $A$ $\in GL(\CC^{N+2})$:
\begin{equation}
\label{map matrix 22}
A=(a_0, ..., a_{N+1})=\begin{bmatrix}
a^{(0)}_0 & a^{(0)}_1 & ... & a^{(0)}_{N+1} \\
a^{(1)}_0 & a^{(1)}_1 & ... & a^{(1)}_{N+1} \\
\vdots & \vdots & \ & \vdots \\
a^{(N+1)}_0 & a^{(N+1)}_1 & ... & a^{(N+1)}_{N+1} \\
\end{bmatrix}\in GL(\CC^{N+2}),\ \
\end{equation}
where each $a_j$ is a column vector in $\CC^{N+2}$, $0\le j\le N+1$. This $A$ is associated to
an automorphism $A^\star$ $\in Aut(\CC\PP^{N+1})$ given by
\begin{equation}
\label{map matrix}
A^\star \bigg(\big[z_0: z_1: ... :z_{N+1}\big]\bigg) =
\bigg[\sum^{N+1}_{j=0} a^{(0)}_j z_j:
\sum^{N+1}_{j=0} a^{(1)}_j z_j: ... :
\sum^{N+1}_{j=0} a^{(N+1)}_j z_j \bigg].
\end{equation}

When $a^{(0)}_0\not=0$, in terms of the non-homogeneous coordinates $(w_1, ..., w_n)$, $A^\star$ is a linear fractional
from $\CC^{N+1}$ which is holomorphic near $(0,...,0)$:
\begin{equation}
\label{map aa}
A^\star\big(w_1, ..., w_{N+1} \big) =
\bigg(
\frac{\sum^{N+1}_{j=0} a^{(1)}_j w_j}{\sum^{N+1}_{j=0} a^{(0)}_jw_j}, ...,
\frac{\sum^{N+1}_{j=0} a^{(N+1)}_j w_j}{\sum^{N+1}_{j=0} a^{(0)}_jw_j}
\bigg),\ \ \ where\ w_j=\frac{z_j}{z_0}.
\end{equation}

We denote $A\in GL^Q(\CC^{N+2})$ if $A$ satisfies $A (Q)\subseteq Q$ where we regard $A$
as a linear transformation of $\CC^{N+2}$. If $A\in GL^Q(\CC^{N+2})$, we must have
$A^\star(\p\HH^{N+1})\subseteq \p\HH^{N+1}$, so that $A^\star\in Aut(\p\HH^{N+1})$.
Conversely, if $A^\star\in Aut(\p\HH^{N+1})$, then $A\in GL^Q(\CC^{N+2})$.


We define a bundle map:
\[
\begin{matrix}
\pi:& GL(\CC^{N+2}) & \ra & \CC\PP^{N+1}\\
\ & A=(a_0, a_1, ..., a_{N+1}) & \mapsto & \pi_0(a_0).
\end{matrix}
\]
Then by (\ref{map matrix}), for any map $A\in GL(\CC^{N+2})$,
$A\in \pi^{-1}\big(\pi_0(a_0)\big)$ $ \Longleftrightarrow$
\ $A^\star ([1:0:...:0])=\pi_0(a_0)$. In particular, by the restriction, we consider a map
\begin{equation}
\label{GLQ pHH}
\begin{matrix}
\pi:& GL^Q(\CC^{N+2}) & \ra & \p\HH^{N+1}\\
\ & A=(a_0, a_1, ..., a_{N+1}) & \mapsto & \pi_0(a_0).
\end{matrix}
\end{equation}
We get $\p \HH^{N+1}\simeq GL^Q(\CC^{N+2})/P_1$ where $P_1$ is the isotropy
subgroup of $GL^Q(\CC^{N+2})$.
Then by (\ref{map matrix}), for any map $A\in GL^Q(\CC^{n+2})$,
\begin{equation}
\label{fiber map}
A\in \pi^{-1}\big(\pi_0(a_0)\big)\ \Longleftrightarrow\
A^\star ([1:0:...:0])=\pi_0(a_0).
\end{equation}

\bigskip

\noindent{\bf CR submanifolds of $\p\HH^{N+1}$}\ \ \ \
Let $H: M' \ra \p \HH^{N+1}$ be a CR smooth embedding
where $M'$ is a strictly pseudoconvex smooth real hypersurface in $\CC^{n+1}$.
We denote $M=H(M')$.

Let $R_{M'}$ be the Reeb vector field of $M'$ with respect to a fixed contact form on $M'$.
Then the real vector $R_{M'}$ generates a real line bundle over $M'$, denoted
by ${\cal R}_{M'}$.
Since we can regard the rank $n$ complex vector bundle $T^{1,0}M'$ as the rank $2n$ real
vector bundle, over the real number field $\RR$ we have:
\begin{equation}
T M'=T^c M'\oplus {\cal R}_{M'} \simeq T^{1,0}M' \oplus {\cal R}_{M'}.
\end{equation}
given by
\begin{equation}
\label{TM basis}
(a_j \frac{\p}{\p x_j}, b_j\frac{\p}{\p y_j})+c R_{M'} \mapsto (a_j +ib_j)\frac{\p }{\p z_j} + c R_{M'}, \ \
\forall a_j, b_j, c\in \RR.
\end{equation}

Since $H$ is a CR embedding, we have
\begin{equation}
H_*(T^{1,0}M')=T^{1,0}M\subset T^{1,0}(\p\HH^{N+1}),
T M \simeq H_*(T^{1,0}M') \oplus H_*({\cal R}_{M'})\subset T(\p\HH^{N+1}).
\end{equation}

\bigskip

\noindent{\bf Lifts of the CR submanifolds}\ \ \ \
Let $M=H(M')\subset \p\HH^{N+1}$ be as above. Consider the commutative  diagram
\[
\begin{matrix}
\ &\ & GL^Q(\CC^{N+2}) \\
\ & e \nearrow  & \downarrow \pi\\
M & \hookrightarrow & \p\HH^{N+1}\\
\end{matrix}
\]
Any map $e$ satisfying $\pi\circ e=Id$ is called a {\it lift} of $M$ to $GL^Q(\CC^{N+2})$.

In order to define a more specific lifts, we need to give some
relationship between geometry on $\p\HH^{N+1}$ and on $\CC^{N+2}$ as follows.
For any subset  $X\in \p\HH^{N+1}$, we denote $\hat X:=\pi_0^{-1}(X)$ where
$\pi_0: \CC^{N+2}-\{0\} \ra \CC\PP^{N+1}$ is the
standard projection map (\ref{pi0}). In particular, for any $x\in M$, $\hat x$ is a complex
line and for the real submanifold $M^{2n+1}$, the real submanifold
$\hat M^{2n+3}$ is of dimension $2n+3$.

For any $x\in M$, we take $v\in \hat x=\pi_0^{-1}(x) \subset \CC^{N+2}-\{0\}$, and we
define
\[
\hat T_x M =T_v \hat M,\  \
\hat T^{1,0}_x M =T^{1,0}_v \hat M, \ \
\hat {\cal R}_{M, x}:={\cal R}_{\hat M, v}
\]
where ${\cal R}_{\hat M}=\cup_{v\in \hat M} {\cal R}_{\hat M, v}$.
These definitions are independent of choice of $v$.

A lift $e=(e_0, e_\alpha, e_\mu, e_{N+1})$ of $M$ into $GL^Q(\CC^{N+2})$, where $1\le \alpha\le n$ and
$n+1\le \mu\le N$, is called a {\it first-order adapted lift}  if it satisfies the conditions:
\begin{equation}
\label{adapted lift}
e_0(x)\in \pi_0^{-1}(x),\ \
span_\CC(e_0, e_\alpha)(x)=\hat T^{1,0}_x M,\ \ span(e_0, e_\alpha, e_{N+1})(x)=\hat T^{1,0}_x M \oplus \hat {\cal R}_{M, x}
\end{equation}
where
\begin{equation}
\label{adapted lift c}
span(e_0, e_\alpha, e_{N+1})(x):=\{c_0 e_0 + c_\alpha e_\alpha+c_{N+1} e_{N+1}\ |\ c_0, c_\alpha\in \CC,\
c_{N+1}\in \RR\}.
\end{equation}
Here we used (\ref{TM basis}) and the fact that the Reeb vector is real.
Locally first-order adapted lifts always exist (see Theorem \ref{existence of lift} below).

\bigskip

We have the restriction bundle ${\cal F}_M^0: =
GL^Q(\CC^{N+2})|_M$ over $M$.
The subbundle $\pi: {\cal F}^1_M \rightarrow M$ of
${\cal F}^0_M$ is defined by
\begin{eqnarray*}
{\cal F}^1_M =\{(e_0, e_j, e_\mu, e_{N+1})\in {\cal F}^0_M\ |\
\ [e_0]\in M,\ (\ref{adapted lift})\ are\ satisfied\}.
\end{eqnarray*}
Local sections of ${\cal F}^1_M$ are exactly all local first-order adapted lifts of $M$.

For two first-order adapted lifts $s=(e_0, e_j, e_\mu, e_{N+1})$
and $\w s=(\w e_0, \w e_j, \w e_\mu, \w e_{N+1})$, by (\ref{adapted lift}), we have
\begin{equation}
\label{e and w e}
\begin{cases}
\w e_0=g^0_0 e_0,\\
\w e_j=g^0_j e_0 + g^k_j e_k,\\
\w e_\mu=g^0_\mu e_0+g^j_\mu e_j + g^\nu_\mu e_\nu + g^{N+1}_\mu e_{N+1},\\
\w e_{N+1}=g^0_{N+1} e_0 + g^j_{N+1} e_j + g^{N+1}_{N+1} e_{N+1},
\end{cases}
\end{equation}
Notice that by (\ref{TM basis}), $g^{N+1}_{N+1}$ is some real-valued function, while other
are complex-valued functions.
In other words, $\w s= s \cdot g$ where
\begin{equation}
\label{g}
g=(g_0, g_j, g_\mu, g_{N+1})=\begin{pmatrix}
g^0_0 & g^0_k & g^0_\mu & g^0_{N+1}\\
0 & g^j_k & g^j_\mu & g^j_{N+1}\\
0 & 0 & g^\nu_\mu & 0\\
0 & 0 & g^{N+1}_\mu & g^{N+1}_{N+1}
\end{pmatrix}
\end{equation}
is a smooth map from $M$ into $GL^Q(\CC^{N+2})$.
Then the fiber of $\pi: {\cal F}^1_M\ra M$ over a point is
isomorphic to the group
\[
G_1 = \bigg\{ g=
\begin{pmatrix}
g^0_0 & g^0_\beta & g^0_\mu & g^0_{N+1}\\
0 & g^\alpha_\beta & g^\alpha_\mu & g^\alpha_{N+1}\\
0 & 0 & g^\nu_\mu & 0\\
0 & 0 & g^{N+1}_\mu & g^{N+1}_{N+1}
\end{pmatrix}\in GL^Q(\CC^{N+2})\bigg\},
\]
where we use the index ranges $1\le \alpha, \beta\le n$ and $n+1\le \mu,
\nu\le N$.

We pull back the Maurer-Cartan form from $GL^Q(\CC^{N+2})$ to ${\cal F}^1_M$ by a
first-order adapted lift $e$ of $M$ as
\[
\omega
=\begin{pmatrix}
\omega^0_0 & \omega^0_\beta & \omega^0_\nu & \omega^0_{N+1}\\
\omega^\alpha_0 & \omega^\alpha_\beta & \omega^\alpha_\nu &
\omega^\alpha_{N+1}
\\
\omega^\mu_0 & \omega^\mu_\beta & \omega^\mu_\nu& \omega^\mu_{N+1}\\
\omega^{N+1}_0 & \omega^{N+1}_\beta & \omega^{N+1}_\nu&
\omega^{N+1}_{N+1}\\
\end{pmatrix}.
\]
Since $\omega=e^{-1}de$, i.e., $e \omega = d e$. Then we have
\begin{equation}
\label{d e0}
d e_0 = e_0 \omega^0_0 + e_\alpha \omega^\alpha_0 + e_\mu \omega^\mu_0
+ e_{N+1}\omega^{N+1}_0.
\end{equation}
On the other hand, we claim:
\begin{equation}
\label{claim d e0}
d e_0=e_0 \omega^0_0 + e_\alpha \omega^\alpha_0 + e_{N+1}\omega^{N+1}_0.
\end{equation}
In fact, take local coordinates systems $(x_1, ..., x_{2n+1})$ for the
real manifold $M$, and $(y_1, y_2, x_1$, $..., x_{2n+1})$ for the
real manifold $\hat M$ where $(y_1, y_2)$ is the coordinates for fibers.
By the first condition in (\ref{adapted lift}),
fixing $x_1, ..., x_{j-1}, x_{j+1}, ...,$ $x_{2n+1}$, $e_0(..., x_j,...)$ is a curve
into $M$ with parameter $x_j$. Then $\frac{\p e_0}{\p x_j} \in T \hat M$ is a tangent
vector to this curve.
Since
$span(e_0, e_\alpha, e_{N+1})(x)=\hat T^{1,0}_x M \oplus \hat {\cal R}_{M, x}$
in (\ref{adapted lift}) and $T \hat M\cong T^{1,0} \hat M \oplus {\cal R}_{\hat M}$,
we obtain
\begin{equation}
\label{e0 xj}
\frac{\p e_0}{\p x_j} = b^j_0 e_0 + b^j_\alpha e_\alpha + b^j_{N+1} e_{N+1}, \ \ 1\le j\le 2n+1
\end{equation}
for some functions $b^j_0, b^j_\alpha$ and $b^j_{N+1}$. We also have
\begin{equation}
\label{e0 yi}
\frac{\p e_0}{\p y_i} = 0,\ \ for\ i = 1, 2,
\end{equation}
because $(y_1, y_2)$ are the coordinates for fibers.
From (\ref{e0 xj}) and (\ref{e0 yi}), we get
\[
d e_0 = \frac{\p e_0}{\p y_1} dy_1 + \frac{\p e_0}{\p y_2} d y_2
+ \sum_j \frac{\p e_0}{\p x_j} d x_j  = \sum_j(b^j_0 e_0 + b^j_\alpha e_\alpha + b^j_{N+1} e_{N+1} )dx_j
\]
\begin{equation}
\label{d e0 00}
=(\sum_j b^j_0 d x_j)e_0+(\sum_j b^j_\alpha d x_j)e_\alpha + (\sum_j b^j_{N+1} d x_j) e_{N+1}.
\end{equation}
Since the $1$-forms $\omega^0_0, \omega^\alpha_0, \omega^\alpha_{N+1}$ in  (\ref{d e0}) are unique, from
(\ref{d e0 00}), it proves Claim (\ref{claim d e0}).

By (\ref{d e0}) and (\ref{claim d e0}), we conclude $\omega^\mu_0=0$, $\forall \mu$.
By the Maurer-Cartan equation
$d\omega=-\omega\wedge \omega$, one gets $0=d\omega^\nu_0=-\omega^\nu_\alpha
\wedge \omega^\alpha_0 - \omega^\nu_{N+1} \wedge \omega^{N+1}_0$, i.e.,
$0 = - \omega^\nu_\alpha \wedge \omega^\alpha_0,\ mod(\omega^{N+1}_0)$.
Then by Cartan's lemma,
\[
\omega^\nu_\beta=q^\nu_{\alpha \beta}\omega^\alpha_0 \ \ mod(\omega^{N+1}_0), \]
for some functions $q^\nu_{\alpha\beta}=q^\nu_{\beta\alpha}$.

\bigskip

\noindent{\bf The CR second fundamental form}\ \ \ \
In order to define the CR second fundamental form $II_M=II^s_M=
q^\mu_{\alpha \beta}\omega^\alpha_0\omega^\beta_0\otimes \underline{e}_\mu$,
mod$(\omega^{N+1}_0)$, let us define $\underline{e}_\mu$ as follows.

For any first-order adapted  lift $e=(e_0, e_\alpha, e_\nu, e_{N+1})$
with $\pi_0(e_0)=x$,  we have
$e_\alpha\in \hat T^{1,0}_x M$.
Recall $T_E G(k, V) \simeq E^*\otimes (V/E)$
where $G(k, V)$ is  the Grassmannian of $k$-planes that pass through the origin in a
vector
space $V$ over $\RR$ or $\CC$ and $E\in G(k, V)$ (\cite{IL03}, p.73). Then $T_x M \simeq (\hat x)^* \otimes
(\hat T_x M/\hat x)$ and hence the vector $e_\alpha$ induces
$\underline{e_\alpha} \in T^{1,0}_x M$
by
\begin{equation*}
\underline{e}_\alpha = e^0\otimes \big(e_\alpha\ mod(e_0) \big),
\end{equation*}
where we denote by $(e^0, e^\alpha, e^\mu, e^{N+1})$ the dual basis of $(\CC^{N+2})^*$.
Similarly, we let
\begin{equation}
\label{underline e mu CR -1}
\underline{e}_\mu=e^0\otimes\big(e_\mu\ mod\ \hat T^{(1,0)}_x M\big) \in
N^{1,0}_x M,
\end{equation}
where $N^{1,0} M$ is the CR normal bundle of $M$ defined by
$N^{1,0}_x M = T^{1,0}_x (\p\HH^{N+1})/ T^{1,0}_x M$.


By direct computation, we obtain a tensor
\begin{equation}
\label{II M 3}
II_M = II_M^e = q^\mu_{\alpha\beta} \omega^\alpha_0 \omega^\beta_0 \otimes
\underline{e}_\mu \in
\Gamma\big(M, S^2 T^{1,0*}_{\pi_0(e_0)} M \otimes N^{1,0}_{\pi_0(e_0)} M\big) \ \
mod(\omega^{N+1}_0).
\end{equation}
The tensor $II_M$ is called the {\it CR second fundamental form} of $M$.

\bigskip

\noindent{\bf Pulling back a lift}\ \ \ \
Let $M\subset \p\HH^{N+1}$ be as above with a point $Q_0\in M$.
Let $A\in GL^Q(\CC^{N+2})$, $A^\star\in Aut(\p\HH^{N+1})$ with $A^\star(Q_0)=P_0$
and $\w M=A^\star(M)$.
Let $\w s: \w M \ra GL^Q(\CC^{N+2})$ be a lift. We claim:
\begin{equation}
\label{pull back lift}
s:=A^{-1} \cdot\w s \circ A^\star,
\end{equation}
is also a lift from $M$ into $GL^Q(\CC^{N+2})$.
In fact, in order to prove that $s$ is a lift, it suffices to
prove: $\pi s = Id$, i.e., for any point $Q\in M$ near $Q_0$,  $\pi s(Q)=Q$.
In fact,
\[
\pi s(Q)=\pi(A^{-1} \cdot \w s \circ A^\star)(Q)=\pi(A^{-1} \cdot \w s(P))=(A^\star)^{-1}(\pi \w s(P))
=(A^\star)^{-1}(P)=Q.
\]
so that our claim is proved.

If, in addition, $\w s$ is a first-order adapted lift of
$\w M$ into $GL^Q(\CC^{N+2})$, $s$ is also a first-order adapted lift of
$M$ into $GL^Q(\CC^{N+2})$.

Let $\Omega$ be the Maurer-Cartan form over
$GL^Q(\CC^{N+2})$. Then by the invariant property $A^* \Omega=\Omega$, we have
$s^* \Omega =(A^{-1} \cdot \w s \circ A^\star)^* \Omega=(A^\star)^* (\w s)^* (A^{-1})^*\Omega
=(A^\star)^*(\w s)^* \Omega$, i.e., it holds on $M$ that
\begin{equation}
\label{omega and w omega}
\omega=(A^\star)^* \w \omega
\end{equation}
where $\omega=s^*\Omega$ and $\w \omega=\w s^*\Omega$ so that
$\omega^\alpha_0=(A^\star)^* \w \omega^\alpha_0$ and $\omega^\mu_\beta=(A^\star)^* \w
\omega^\mu_\beta$.
The last equality yields
\begin{equation}
\label{s and w s}
q^\mu_{\alpha \beta} = \w q^\mu_{\alpha \beta}\circ A^\star.
\end{equation}

\bigskip

\section{CR second fundamental forms ----- Definition 4}

\

Definition 4 will be the one as a tensor with respect to the group $SU(N+1,1)$.

As for Definition 3, we consider the real hypersurface $Q$ in $\CC^{N+2}$ defined
by the homogeneous equation
\begin{equation}
\label{3 sphere}
\langle Z, Z\rangle : = \sum_A Z^A \ov{Z^A} + \frac{i}{2}(Z^{N+1}\ov{Z^0}
- Z^0
\ov{Z^{N+1}})=0,
\end{equation}
where $Z=(Z^0, Z^A, Z^{N+1})^t\in \CC^{N+2}$. This can be extended
to the scalar product
\begin{equation}
\label{product}
\langle Z, Z'\rangle : = \sum_A Z^A \ov{{Z'}^A} +
\frac{i}{2}(Z^{N+1} \ov{Z'}^0
- Z^0\ov{{Z'}^{N+1}}),
\end{equation}
for any $Z=(Z^0, Z^A, Z^{N+1})^t, Z'=({Z'}^0, {Z'}^A, {Z'}^{N+1})^t \in \CC^{N+2}$.
This product has the properties:
$\langle Z, Z'\rangle$ is linear in $Z$ and anti-linear in $Z'$; $\ov{\langle Z, Z'\rangle}
=\langle Z', Z\rangle$; and $Q$ is defined by $\langle Z, Z\rangle=0$.

Let $SU(N+1, 1)$ be the group of unimodular linear transformations
of $\CC^{N+2}$ that leave the form $\langle Z, Z\rangle$ invariant
(cf. [CM74]).

By a {\it $Q$-frame} is meant an element $E=(E_0, E_A$, $E_{N+1})\in
GL(\CC^{N+2})$ satisfying (cf. [CM74, (1.10)])
\begin{equation}
\label{3 basic formula theta}
\left\{
\begin{array}{lll}
&& det(E)=1,\\
&&\langle E_A, E_B\rangle = \delta_{AB},\
\langle E_0, E_{N+1}\rangle=-\langle E_{N+1},  E_0
\rangle = - \frac{i}{2},\\
\end{array}
\right.
\end{equation}
while all other products are zero.

There is exactly one transformation of $SU(N+1, 1)$ which maps a
given $Q$-frame into another.
By fixing one $Q$-frame as reference, the group $SU(N+1, 1)$ can be identified with the
space of all $Q$-frames. Then $SU(N+1,1)\subset GL^Q(\CC^{N+1})$ is a subgroup with the composition
operation. By  (\ref{GLQ pHH}) and the restriction, we have the projection
\begin{equation}
\pi: SU(N+1, 1) \ra \p\HH^{N+1},\ (Z_0, Z_A, Z_{N+1})  \mapsto  span(Z_0).
\end{equation}
which is called a {\it $Q$-frames bundle}. We get $\p \HH^{N+1}\simeq SU(N+1, 1)/P_2$ where $P_2$ is the isotropy
subgroup of $SU(N+1, 1)$. $SU(N+1, 1)$ acts on $\p\HH^{N+1}$ effectively.

\bigskip

Consider $E=(E_0, E_A, E_{N+1})\in SU(N+1, 1)$ as a local lift.
Then the {\it Maurer-Cartan form} $\Theta$ on $SU(N+1, 1)$ is defined by
$d E = (d E_0, d E_A, d E_{N+1}) = E\Theta$, or  $\Theta= E^{-1} \cdot d
E$,
i.e.,
\begin{equation}
\label{3 structure eq}
d\begin{pmatrix}
E_0 & E_A & E_{N+1}\end{pmatrix}
= \begin{pmatrix}
E_0 & E_B & E_{N+1}
\end{pmatrix}
\begin{pmatrix}
\Theta^0_0 & \Theta^0_A & \Theta^0_{N+1}\\
\Theta^B_0 & \Theta^B_A & \Theta^B_{N+1}\\
\Theta^{N+1}_0 & \Theta^{N+1}_A & \Theta^{N+1}_{N+1}\\
\end{pmatrix},
\end{equation}
where $\Theta^B_A$ are 1-forms on $SU(N+1, 1)$.
By (\ref{3 basic formula theta}) and (\ref{3 structure eq}),
the Maurer-Cartan form $(\Theta)$ satisfies
\begin{equation}
\label{3 basic formula theta 2}
\begin{array}{lll}
&& \Theta^0_0 + \ov{\Theta^{N+1}_{N+1}} =0,\
\Theta^{N+1}_0=\ov{{\Theta}_0^{N+1}},\
\Theta^0_{N+1}=\ov{\Theta_{N+1}^0},\\
&&
\Theta^{N+1}_A=  2i \ov{\Theta^A_0},\ \Theta^A_{N+1}=- \frac{i}{2}
\ov{\Theta^0_A},\
\Theta^A_B+\ov{\Theta^B_A}=0,\ \Theta^0_0+\Theta^A_A+\Theta^{N+1}_{N+1}= 0,
\end{array}
\end{equation}
where $1\le A\le N$. For example, from $\langle E_A, E_B\rangle=\delta_{A B}$, by
taking differentiation, we obtain
\[
\langle dE_A, E_B\rangle+\langle E_A, dE_B\rangle=0.
\]
By (\ref{3 structure eq}), we have
\[
\begin{cases}
d E_0=E_0\Theta^0_0+E_B \Theta^B_0+E_{N+1}\Theta^{N+1}_0,\\
d E_A=E_0\Theta^0_A+E_B \Theta^B_A+E_{N+1}\Theta^{N+1}_A,\\
d E_{N+1}=E_0\Theta^0_{N+1}+E_B \Theta^B_{N+1}+E_{N+1}\Theta^{N+1}_{N+1}.\\
\end{cases}
\]
Then
\[
\langle E_0 \Theta^0_A+E_C \Theta^C_A + E_{N+1}\Theta^{N+1}_A,\ E_B\rangle
+\langle E_A, E_0 \Theta^0_B + E_D \Theta^D_B + E_{N+1} \Theta^{N+1}_B\rangle=0,
\]
which implies $\Theta^B_A+\ov{\Theta^A_B}=0$. In particular,
from (\ref{3 basic formula theta 2}), $\Theta^0_A=- 2i
\ov{\Theta^A_{N+1}}$.
$\Theta$ satisfies
\begin{equation}
\label{15}
d\Theta=-\Theta\wedge \Theta.
\end{equation}

Let $M \hookrightarrow \p\HH^{N+1}$ be the image of $H:M'\ra
\p\HH^{N+1}$ where $M'\subset \CC^{n+1}$ is a CR strictly
pseudoconvex smooth hypersurface. Consider the
inclusion map $M \hookrightarrow \p\HH^{N+1}$ and a lift $e=(e_0,
e_1, ..., e_{N+1})=(e_0, e_\alpha, e_\nu$, $e_{N+1})$ of $M$ where
$1\le \alpha\le n$ and $n+1\le \nu \le N$
\[
\begin{matrix}
\ &\ & SU(N+1, 1)\\
\ & e \nearrow & \downarrow \pi \\
M & \hookrightarrow & \p\HH^{N+1}\\
\end{matrix}
\]
We call $e$ a {\it first-order adapted lift} if for any $x\in M$,
\begin{equation}
\label{adapted lift 2}
\pi_0\big(e_0(x)\big)=x,\
span_\CC(e_0, e_\alpha)(x)=\hat T^{1,0}_x M,\ span(e_0, e_\alpha, e_{N+1})(x)
=\hat T^{1,0}_x M \oplus \hat{\cal R}_{M, x}.
\end{equation}
Locally first-order adapted lifts always exist (see Theorem \ref{existence of lift} below).
We have the restriction bundle ${\cal F}_M^0: =
SU(N+1,1)|_M$ over $M$.
The subbundle $\pi: {\cal F}^1_M \rightarrow M$ of
${\cal F}^0_M$ is defined by
\begin{eqnarray*}
{\cal F}^1_M =\{(e_0, e_j, e_\mu, e_{N+1})\in {\cal F}^0_M\ |\
\ [e_0]\in M,\ (\ref{adapted lift 2})\ are\ satisfied\}.
\end{eqnarray*}
Local sections of ${\cal F}^1_M$ are exactly all local first-order adapted lifts of $M$.
The fiber of $\pi: {\cal F}^1_M\ra M$ over a point is
isomorphic to the group
\[
G_1 = \bigg\{ g=
\begin{pmatrix}
g^0_0 & g^0_\beta & g^0_\nu & g^0_{N+1}\\
0 & g^\alpha_\beta & g^\alpha_\nu & g^\alpha_{N+1}\\
0 & 0 & g^\mu_\nu & 0\\
0 & 0 & 0 & g^{N+1}_{N+1}
\end{pmatrix}\in SU(N+1, 1)\bigg\},
\]
where we use the index ranges $1\le \alpha, \beta\le n$ and $n+1\le
\mu, \nu \le N$.


By the remark below (\ref{e and w e}), $g_{N+1}^{N+1}$ is real-valued.
By (\ref{3 basic formula theta}), we have $\langle
g_0, g_{N+1}\rangle=-\frac{i}{2}$, it implies
$g^0_0 \cdot \ov{g^{N+1}_{N+1}}=1$. In particular, both $g^{N+1}_{N+1}$ and $g^0_0$ are real.
Since $\langle g_0, g_\mu\rangle=0$ and $g^0_0\not=0$, it implies
$g^{N+1}_\mu=0$. Since $\langle g_\alpha, g_\beta \rangle = \delta_{\alpha\beta}$, it implies that the matrix
$(g^\beta_\alpha)$ is unitary. Since $deg(g)=1$, it implies $g^0_0 \cdot det(g^\beta_\alpha)\cdot
det(g^\nu_\mu)\cdot g^{N+1}_{N+1}=1$, i.e.,
$det(g^\beta_\alpha)\cdot det(g^\nu_\mu)=1$.

\bigskip

By considering all first-order adapted lifts from $M$ into $SU(N+1,1)$, as the definition of $II_M$
in Definition 3, we can defined CR second fundamental form $II_M$ as
in (\ref{II M 3}):
\begin{equation}
\label{II M 4 2}
II_M = II^e_M = q^\mu_{\alpha\beta} \omega^\alpha_0 \omega^\beta_0 \otimes
\underline{e}_\mu \in
\Gamma(M, S^2 T^{1,0*}_{\pi_0(e_0)} M \otimes N^{1,0}_{\pi_0(e_0)} M),\ \
mod(\omega^{N+1}_0),
\end{equation}
which is a well-defined tensor, and is called the {\it CR second
fundamental form} of $M$.

We remark that the notion of $II_M$ in Definition 4 was introduced in a paper by S.H.
Wang \cite{Wa06}.

\bigskip


\noindent{\bf Pulling back a lift}\ \ \ \
Let $M\subset \p\HH^{N+1}$ be as above with a point $Q_0\in M$.
Let $A\in SU(N+1,1)$, $A^\star\in Aut(\p\HH^{N+1})$ with $A^\star(Q_0)=P_0$
and $\w M=A^\star(M)$.
Let $\w s: \w M \ra SU(N+1,1)$ be a lift. We claim:
\begin{equation}
\label{pull back lift, 4}
s:=A^{-1} \cdot\w s \circ A^\star,
\end{equation}
is also a lift from $M$ into $SU(N+1,1)$.
Similarly as in (\ref{omega and w omega}) and (\ref{s and w s}), we have

\begin{equation}
\label{omega and w omega, 4}
\omega=(A^\star)^* \w \omega
\end{equation}
and
\begin{equation}
\label{s and w s, 4}
q^\mu_{\alpha \beta} = \w q^\mu_{\alpha \beta}\circ A^\star.
\end{equation}
where $\omega=s^*\Omega$, $\w \omega=\w s^*\Omega$ and $\Omega$ is the Maurer-Cartan
form over $SU(N+1,1)$.


\noindent{\bf [Example]}\ \ \ \ Consider the maps in (\ref{sigma}) and
(\ref{tau}):
\begin{eqnarray*}
&&\sigma^0_p(z,w)=(z+z_0, w+w_0+2i \langle z,\overline{z_0}
\rangle),\\
&&\tau^F_p(z^*,w^*)=(z^*-\widetilde{f}(z_0,w_0),w^*-\overline{g(z_0,w_0)}-
2i \langle z^*,\overline{\widetilde{f}(z_0,w_0)} \rangle )
\end{eqnarray*}
where $p=(z_0, w_0)$, $z=\CC^n$, $w=z_{n+1}$, $\sigma^0_p\in
Aut(\p\HH^{n+1})$, and
$\tau^F_p\in Aut(\p\HH^{N+1})$.


By (\ref{map matrix 22}) and (\ref{map aa}), these two maps correspond to
two matrices:
\begin{equation}
\label{ex 1}
A_{\sigma^0_p}=\begin{bmatrix}
1 & 0 & ... & 0 & 0 \\
z_{0 1} & 1 & ... & 0 & 0 \\
\vdots & \vdots & \ddots & \vdots & \vdots \\
z_{0 n} & 0 & ... & 1 & 0 \\
w_0 & 2i \ov{z_{01}} & ... & 2i \ov{z_{0 n}} & 1 \\
\end{bmatrix}\in SU(n+1,1)
\end{equation}
and
\begin{equation}
\label{ex 1a}
A_{\sigma^F_p}=
\begin{bmatrix}
1 & 0 & ... & 0 & 0 \\
-\w f_{0 1} & 1 & ... & 0 & 0 \\
\vdots & \vdots & \ddots & \vdots & \vdots \\
-\w f_{0 N-n} & 0 & ... & 1 & 0 \\
- \ov{g(z_0, w)} & -2i \ov{\w f_1(z_0, w_0)} & ... & -2i \ov{\w
f_{N-n}(z_0, w_0)} & 1 \\
\end{bmatrix} \in SU(N+1,1)
\end{equation}
where $z_0=(z_{01}, ..., z_{0n})$ and $w_0=z_{0 n+1}$.\ \ \ $\Box$
\bigskip


\noindent{\bf [Example]}\ \ \ \ Consider the map $F_{\lambda, r, \vec a,
U}=(f, g)\in
Aut_0(\p\HH^{n+1})$
\[
f(z)=\frac{\lambda(z+\vec aw)U}{1-2i\langle z, \ov{\vec
a}\rangle-(r+i\Vert \vec
a\Vert^2)w},\
g(z)=\frac{\lambda^2 w }{1-2i\langle z, \ov{\vec a}\rangle-(r+i\Vert \vec
a\Vert^2)w}
\]
where $\lambda>0, r\in \RR, \vec a\in \CC^n$ and $U=(u_{\alpha \beta})$ is
an
$(n-1)\times(n-1)$ unitary matrix.
By (\ref{map matrix 22}) and (\ref{map aa}), its corresponding matrix,
\begin{equation}
\label{ex 2}
A_{F_{\lambda, r, \vec a, U}}
=
\begin{bmatrix}
1 & -2i \ov{a_1} & ... & -2i\ov{a_n} & -(r+i\Vert \vec a\Vert^2) \\
0 & \lambda u_{11} & ... & \lambda u_{1 n}& \lambda a_1\\
\vdots & \vdots & \ddots & \vdots & \vdots \\
0 & \lambda u_{n 1} & ... & \lambda u_{n n} & \lambda a_n \\
0 & 0 & ... & 0 & \lambda^2\\
\end{bmatrix},
\end{equation}
is not in $SU(n+1,1)$ in general. In fact, we can write
\begin{equation}
\label{6 composition}
F_{\lambda, r, \vec a, U}=F_{\lambda, 0, 0, Id} \circ F_{1, 0, 0, U} \circ
F_{1, r, \vec
a, Id}.
\end{equation}
or
$
A_{F_{\lambda, r, \vec a, U}} = A_{F_{\lambda, 0, 0, Id}} \cdot A_{F_{1,
0, 0, U}} \cdot
A_{ F_{1, r, \vec a, Id}}.
$
 Here $A_{F_{1, 0, 0, U}}$ and $A_{F_{1, r, \vec a, Id}}$ are in
$SU(N+1,1)$; while
$A_{F_{\lambda, 0, 0, Id}}$  is in
$SU(N+1,1)$ if and only if $\lambda =1$.
Therefore
\begin{equation}
\label{7, example}
A_{F_{\lambda, r, \vec a, U}}\ is\ in\ SU(n+1,1)\ if\ and\ only\ if\
\lambda=1.
\end{equation}

\bigskip

\section{Existence of First-order Adapted Lifts from $M$ into
$SU(N+1,1)$ or into  $GL^Q(\CC^{N+2})$}


\noindent{\bf Existence of first-order adapted lifts.} \ \ \ \
 Let
$(M', 0)$  be a germ of smooth real hypersurface in $\CC^{n+1}$ defined by the defining function
\begin{equation}
r=\sum^n_{j=1} z_j \ov z_j + \frac{i}{2}(w - \ov w)+o(2).
\end{equation}
We take
\[
\theta = i \p r=i \bigg( \sum^{n}_{j=1} \ov{z_j} d z_j - \frac{1}{2} d w)\bigg) + o(1).
\]
as a contact form of $M'$.

Write $w=u+i v$. Here $v=\sum^n_{j=1} |z_j|^2+o(2)$. Take $(z_j, u)$ as a coordinates system of $M'$.
By considering the coordinate map: $h: \CC^n\times \RR \ra M',\ (z_j, u) \mapsto (z_j, u+i|z|^2+o(2))$,
we get the pushforward
\[
h_*(\frac{\p}{\p z_j})=L_j:=\frac{\p }{\p z_j} + i \big(\ov{z_j} + o(1)\big)\frac{\p }{\p u},\ \
h_*(\frac{\p}{\p u})=R_{M'}:=(1+o(1)) \frac{\p }{\p u}
\]
for $j=1,2,...,n$. Then $\{L_j\}_{1\le j\le n}$ form a basis of the complex tangent bundle
$T^{1,0} M'$ of $M'$. Since $d \alpha = - i \sum^n_{j=1} dz_j \wedge d \ov{z_j}$,  we see
that $R$ is the Reeb vector field of $M'$.  In particular, as the restriction at $0$, we have
\begin{equation}
\label{Lj R at 0}
L_j|_0=\frac{\p}{\p z_j}|_0,\ \ R_{M'}|_0=\frac{\p }{\p u}|_0.
\end{equation}

\bigskip

\begin{theorem}
\label{existence of lift}
Let $M \hookrightarrow \p\HH^{N+1}$ be the image of $H:M'\ra
\p\HH^{N+1}$ where $M'\subset \CC^{n+1}$ is a smooth strictly
pseudoconvex CR-hypersurface. Then for any point in $M$, the
first-order adapted lift $E=(E_0, E_\alpha, E_\mu, E_{N+1})$ of $M$ into $SU(N+1,1)$ ( hence into $GL^Q(\CC^{N+2})$)
exists in some neighborhood of the point in $M$.
\end{theorem}

\bigskip

\noindent{\it Proof:}\ \ {\bf Step 1.}\ \ \ Without of loss of generality, we
assume that $0\in M$ so that it suffices to construct a lift $E=(E_0, E_\alpha, E_\mu,
E_{N+1})$ in a neighborhood of the point $0$. Here we denote $[1:0:...:0]$ by 0.

Assume that $M'$ is defined by the equation $Im\ w = |z|^2 + o(|z|^2)$ in
$(z, w)\in \CC^n \times \CC$ where $w=u+iv$.
Assume that  $H=(1, f_\alpha, \phi_\mu, g)$ is the smooth CR embedding
of $M'$ into $\p\HH^{N+1}$ with $H(0)=0$ and
\begin{equation}
\label{7 normal}
f=z+O(|(z,w)|^2), \phi=O(|(z,w)|^2),\ g=w+O(|(z,w)|^2).
\end{equation}
Let $L_\alpha, \alpha=1,2,...,n$ be a basis of the CR vector fields and $R$ is the Reeb
vector field on $M'$. Then as in (\ref{Lj R at 0}) with (\ref{7 normal}), we have
\begin{equation}
\label{7 at 0}
L_\alpha|_0=\frac{\p}{\p z_j}|_0,\ \  and \ \ R|_0=\frac{\p }{\p u}|_0.
\end{equation}
It follows that $\bar L_\alpha
H=0$ as $H$ is a CR map. By the Lewy extension theorem, $H$ extends
holomorphically to one side of $M'$, denoted by $D$, where $D$ is
obtained by attaching the holomorphic discs. By applying the maximum
principle and the Hopf lemma to the subharmonic function $\sum
|f_\alpha|^2 + \sum |\phi_\mu|^2 + \frac{i}{2}(g-\bar g)$ on $D$,
it follows that $\frac{\partial Im\ g}{\partial v}(0) \neq 0$. Since
$\frac{\partial g}{\partial \bar w}=0$ and $\frac{\partial Im\ g}{\partial u}(0)
=0$, we have $R g(0)=\frac{\partial g}{\partial u}(0)=
\frac{\partial Im\ g}{\partial v}(0) \neq 0$.

\bigskip

\noindent{\bf Step 2. Direct construction of $E_0, E_\alpha$ and $E_{N+1}$}\ \ \ \
We define
\begin{equation}
\label{E0}
E_0:=\begin{bmatrix}
1 \\
f_\alpha(z, w)\\
\phi_\mu(z, w)\\
g(z, w)
\end{bmatrix}
\end{equation}
which can be regarded as a point in $\p\HH^{N+1}$.  Then $\langle E_0, E_0\rangle=0$ holds:
\begin{equation}
\label{E0 E0 0}
\sum f_\alpha \bar f_\alpha + \sum \phi_\mu \bar \phi_\mu +
\frac{i}{2} (g - \bar g) =0,\ \ \ on\ M.
\end{equation}

Apply the CR vector field $L_\beta$ to $E_0$, we define
\[ \w E_\beta= (0, L_\beta f_\alpha, L_\beta \phi_\mu,
L_\beta g)^t, \]
which form the basis of the complex tangent bundle $T^{1,0}_{\pi_0(E_0)}(M)$.
Then in a neighborhood of $0$ in $M$, we have
\[
span_\CC(E_0, \w E_\alpha)=\hat T^{(1,0)}_{\pi_0(E_0)} M.
\]

Now, we have $\langle E_0, \w E_\alpha\rangle =0$ by
applying $L_\beta$ to (\ref{E0 E0 0}):
\begin{equation}
\label{Ea E0}
\sum \bar f_\alpha L_\beta f_\alpha + \sum \bar \phi_\mu L_\beta \phi_\mu +
\frac{i}{2} L_\beta g =0 .
\end{equation}

By the Gram-Schmid orthonormalization procedure, we can obtain,
from $\{\w E_\beta \}$,  an orthonormal set with respect to the usual Hermitian
inner product $\langle\ ,\ \rangle_0$; we denote it by $\{E_\beta \}$.
By the definition (\ref{product}), we notice that for any $Z=(Z^0, Z^A, Z^{N+1})$ and
$ Z'=(Z'^0, Z'^A, Z'^{N+1})$,
\begin{equation}
\label{< >0}
\langle Z,\ Z'\rangle=\bigg\langle (\frac{i}{2} Z^{N+1}, Z^A, -\frac{i}{2} Z^0),\
(Z_0', Z'^A, Z'^{N+1}) \bigg\rangle_0=\langle \hat Z, Z'\rangle_0,
\end{equation}
where $\langle\ ,\ \rangle_0$ is the usual Hermitian inner product and
$\hat Z:=(\frac{i}{2} Z^{N+1}, Z^A, -\frac{i}{2} Z^0)$.
Then we see from (\ref{Ea E0}) that
\[
\langle E_0,\ E_\beta\rangle=\bigg\langle (\frac{i}{2}g, f_\alpha, \phi_\mu, -\frac{i}{2}), \ (0, L_\beta f_\alpha,
L_\beta \phi_\mu, L_\beta g)\bigg\rangle_0=0.
\]
Also we observe $\langle E_\alpha, \ E_\beta\rangle=\langle E_\alpha,
\ E_\beta\rangle_0=\delta_{\alpha\beta}$.
Then $\langle E_0, E_0\rangle = 0, \langle E_0, E_\beta\rangle =0$ and
$\langle E_\alpha, E_\beta\rangle = \delta_{\alpha \beta}$ hold.

Applying the Reeb vector field $R$, we define another vector

\[\w E_{N+1}:= (0, R\ f_\alpha, R\ \phi_\mu,
R\  g)^t\] over a neighborhood of $0$ in $M$ such that
\[
span(E_0, E_\alpha, \w E_{N+1})=\hat T_{\pi_0(E_0)} M.
\]

We want to construct
\[
E_{N+1}=A E_0 + B_\alpha E_\alpha + C \w E_{N+1}\]
such that
\[
\langle E_{N+1},\ E_0\rangle=\frac{i}{2}, \
\langle E_\alpha, E_{N+1}\rangle=0,\ and\  \langle E_{N+1}, \ E_{N+1}\rangle=0.
\]

From $\langle E_{N+1}, E_0\rangle=\frac{i}{2}$, we get $\langle A E_0 + B_\alpha E_\alpha + C \w E_{N+1},
\  E_0\rangle=\frac{i}{2}$ so that
\begin{equation}
\label{C}
C=\frac{i}{2 \langle \w E_{N+1}, E_0\rangle}.
\end{equation}
By (\ref{7 at 0}), we notice that
\[
\langle \w E_{N+1}, E_0\rangle|_0 =
\sum \frac{\p f_\alpha}{\p u}|_0 \bar f_\alpha(0)
 + \sum \frac{\p \phi_\mu}{\p u}|_0 \bar\phi_\mu(0) + \frac{i}{2}\frac{\p g}{\p u}|_0
\]
and therefore $\langle \w E_{N+1}, E_0\rangle(0)=\frac{i}{2}R\ g(0) \neq 0$.

From  $\langle E_{N+1}, E_\alpha \rangle=0$, we get $\langle A E_0 + B_\beta E_\beta + C \w E_{N+1},
\  E_\alpha\rangle=0$ so that
\begin{equation}
\label{B alpha}
B_\alpha=- C \delta_{\beta\alpha} \langle \w E_{N+1}, E_\beta \rangle= - C \langle \w E_{N+1}, E_\alpha\rangle.
\end{equation}

From $\langle E_{N+1}, E_{N+1} \rangle=0$, we get $\langle A E_0 + B_\beta E_\beta + C \w E_{N+1},
\ A E_0 + B_\beta E_\beta + C \w E_{N+1}\rangle=0$. Since $C\langle \w E_{N+1}, E_0\rangle = \frac{i}{2}$,
$\ov C\langle E_0, \w E_{N+1}\rangle = -\frac{i}{2}$, $B_\alpha=-C\langle \w E_{N+1}, E_\alpha \rangle$ and
$\ov{B_\alpha}=-\ov C\langle E_\alpha, \w E_{N+1} \rangle$ by (\ref{C}) and (\ref{B alpha}), we obtain
\[
-\frac{i}{2} A + \frac{i}{2} \ov A - \sum_\alpha|B_\alpha|^2 + |C|^2 \langle E_{N+1}, E_{N+1} \rangle=0,
\]
so that
\begin{equation}
\label{Im A}
Im(A)=\sum_\alpha|B_\alpha|^2-|C|^2\langle E_{N+1}, E_{N+1} \rangle.
\end{equation}
Therefore $E_{N+1}$ is determined.

So far we have
$\langle E_0, E_0\rangle = \langle E_{N+1}, E_{N+1}\rangle =\langle E_0, E_\beta\rangle = \langle E_{N+1}, E_\beta\rangle=0$,
$\langle E_\alpha, E_\beta\rangle = \delta_{\alpha \beta}$ and $\langle E_0, E_{N+1}\rangle = -\frac{i}{2}$ hold.

\bigskip

\noindent{\bf Step 3. Construction of $E$}\ \ \ \
From Step 2, at the point $0$, we have vectors
\begin{equation}
\label{E0 Ea E}
E_0|_0=[1:0:...:0], \ E_1|_0=[0:1:0:...:0], ..., E_n|_0=[0:0:...:1: 0:...:0],
\end{equation}
and
\begin{equation}
\label{E0 Ea E2}
\ E_{N+1}|_0=[0:0:...:0:1].
\end{equation}
Therefore we can define $E$ at the point $0$ by
\begin{equation}
E(0):=Id\in SU(N+1,1).
\end{equation}

For any other point $P$ in a small neighborhood of $0$ in $M$,
we are going to define $E(P)\in SU(N+1,1)$ as follows.

Write $H(p)=P$ for some $p\in M'$. Then we take a map $\Psi_P\in SU(N+1,1)$ such that
\[
\Psi_P^\star(P)=0,\ \ T^{1,0}_0 \Psi(M)=span_\CC (E_0|_0, E_\alpha|_0),\ \ and\ \
T_0 \Psi(M)=span (E_0|_0, E_\alpha|_0, E_{N+1}|_0).
\]
where $E_0|_0, E_\alpha|_0$ and $E_{N+1}|_0$ are defined in (\ref{E0 Ea E}) and (\ref{E0 Ea E2}).
The map $\Psi_P$ can be defined as $A_{F_{1, r, \vec a, U}}\circ A_{\sigma^F_p}$ where $A_{\sigma^F_p}
\in SU(N+1,1)$ as in (\ref{ex 1a}) and $A_{F_{1, r, \vec a, U}}\in SU(N+1,1)$ as in (\ref{ex 2}).
Notice in the construction of the normalization $F^{**}$ and $F^{***}$, we can always choose $\lambda=1$
so that (\ref{6 composition}) can be used. $\Psi_P$ is smooth as $P$ varies. Then we define
\begin{equation}
E(P):= (\Psi^\star_P)^* E(0)=(\Psi_P)^{-1} E(0).
\end{equation}
This definition is the same as in (\ref{pull back lift, 4}).
Since $\Psi_P$ is invariant for the Hermitian scalar product
$\langle\ ,\ \rangle$ defined in (\ref{product}) and $E(0)$ satisfies the identities
(\ref{3 basic formula theta}), it implies
that $E(P)$ satisfies the identities (\ref{3 basic formula theta}),
i.e., $E(p)\in SU(N+1,1)$.

As a matrix, we denote $E(P)=(\hat E_0, \hat E_\alpha, \hat E_\mu, \hat E_{N+1})$.
Since the map $\Psi_P$ preserves the CR structures and the tangent vector spaces of $M$ and $\Psi_P(M)$, we have
\[
span_\CC(\hat E_0, \hat E_\alpha)=span_\CC(E_0, E_\alpha)|_P,\ \
span(\hat E_0, \hat E_\alpha, \hat E_{N+1})=span(E_0, E_\alpha, E_{N+1})|_P.
\]
where $E_0, E_\alpha$ and $E_{N+1}$ are constructed in Step 2. We remark that
we can replace $(\hat E_0, \hat E_\alpha$, $\hat E_{N+1})$ by $(E_0, E_\alpha, E_{N+1})$.\ \ \ $\Box$

\bigskip
\noindent{\bf Existence of a more special first-order adapted lifts when $M$ is spherical} \ \ \ \
When $M=F(\p\HH^{n+1})$ where $F\in Prop_2(\HH^{n+1}, \HH^{N+1})$, we can construct a more special first-order
adapted lift of $M$ into $SU(N+1,1)$ as follows (cf. \cite{HJY09}).

Let $F=(f, \phi, g)\in Prop_2(\p\HH^{n+1}, \p\HH^{N+1})$ be any map with $F=F_p^{***}$.
Then $F(0)=0$.
We introduce a local biholomorphic map near the origin
\[F_{f g}:=(f, g): \CC^{n+1} \ra \CC^{n+1},\ (z, z_{N+1}) \mapsto
(f, g)=(\hat z, \hat z_{N+1})\]
with its inverse
\[
F^{-1}_{f g}: \CC^{n+1} \ra \CC^{n+1}, \ (\hat z, \hat z_{N+1}) \mapsto (
(F^{-1}_{f g})^{(1)}, ..., (F_{f g}^{-1})^{(n)}, (F_{f g}^{-1})^{(N+1)})=(z, z_{N+1}).\]

Here we use $(\hat z, \hat z_{N+1})$ as a coordinates system of $M=F(\p\HH^{n+1})$ near $F(0)=0$.
Denote $Proj_{fg}: \CC^{N+1}\ra \CC^{n+1}, (\hat z, \hat z_\mu, \hat z_{N+1}) \mapsto
(\hat z, \hat z_{N+1})$. Then we have $Proj_{fg}\circ F = F_{fg}$:
\[
\begin{matrix}
F: \p\HH^{n+1} & \rightarrow & M \\
\ & \searrow F_{fg} & \downarrow Proj_{fg}

\\
\ &\ & \CC^{n+1}
\end{matrix}
\]
We also have a pair of inverse maps
$F: \p\HH^{n+1} \ra M\ and\ (F_{f g}^{-1})\circ Proj_{fg}: M \ra \p \HH^{n+1}.$

Locally we can regard $M$ as a graph: $F\circ F_{f g}^{-1}: \CC^{n+1}\ra M\subset \CC^{N+2}$:
\[
(\hat z, \hat z_{N+1}) \mapsto \big(\hat z,\ \phi\big((F_{fg})^{-1}( \hat z, \hat
z_{N+1})\big),\ \ \hat z_{N+1}\big)
\]

Now let us define a lift of $M$ into $SU(N+1,1)$
\begin{equation}
\label{CR lift}
e=(e_0, e_\alpha, e_\mu, e_{N+1})\in SU(N+1, 1),\ \ 1\le \alpha\le n, \ n+1\le \mu\le N
\end{equation}
as follows.

We define $e_0: M\hookrightarrow  \CC^{N+2}$ be the inclusion:
\begin{eqnarray}
\label{Ball e0}
e_0(\hat z, \hat z_{N+1})=F\circ F_{f g}^{-1}(\hat z, \hat z_{N+1})=
\bigg[1 : \  \hat z:\ \phi\big((F_{fg})^{-1}( \hat z, \hat z_{N+1})\big):\ \hat z_{N+1}\bigg]^t
\end{eqnarray}
$\forall (\hat z, \hat z_{N+1}) \in \CC^{n+1}$.
We define $e_\alpha: M \ra \CC^{N+2}$ for $1\le \alpha\le n$:
\begin{equation}
\label{Ball ealpha}
e_\alpha := \frac{1}{\sqrt{|L_\alpha f|^2 + |L_\alpha \phi|^2}}\big[0: L_\alpha f: L_\alpha
\phi: L_\alpha g\big]^t \circ F_{f g}^{-1}.
\end{equation}
where $L_\alpha=\frac{\p }{\p z^\alpha}+2i \bar{z}^\alpha \frac{\p }{\p z^{N+1}}$.
By the definition (\ref{product}), we have $\langle e_0, e_0 \rangle=0$
because
$f \cdot \ov{f} + \phi \cdot \ov{\phi} - \frac{1}{2i}(g - \ov g) = \hat z \cdot \ov{\hat z} + \phi\big((F_{fg})^{-1}( \hat z, \hat z_{N+1})\big) \ov{\phi\big((F_{fg})^{-1}( \hat z, \hat z_{N+1})\big)}
 + \frac{i}{2}(\hat z_{N+1} - \ov{\hat z_{N+1}})=0$ holds on $\p\HH^{n+1}$, and
$\langle e_0, e_\alpha\rangle=0$ because $L_\alpha f\cdot \ov
{f}+L_\alpha \phi \cdot \ov{\phi} + \frac{i}{2} L_\alpha g=0$ holds on
$\p\HH^{n+1}$,
and $\langle e_\alpha, e_\beta\rangle=\delta_{\alpha \beta}$
because
$L_\alpha f \cdot \ov{L_\beta f}+L_\alpha \phi\cdot \ov{L_\beta \phi}=0$ holds on
$\p\HH^{n+1}$ for $\alpha\not=\beta$.

If we define $\w e_{N+1}:=(0, Tf, T\phi, T g)^t\circ F_{f g}^{-1}$, where $T=\frac{\p}{\p u}$
with $z^{N+1}=u+iv$, then $span(e_0, e_\alpha, \ov{e_\alpha}, \w e_{N+1})= \hat T_{\pi_0(e_0)}
M$. We then find coefficient functions $A, B_\alpha$ and $C$
such that $e_{N+1}=A e_0 + \sum B_\alpha e_\alpha + C \w e_{N+1}$ satisfies
\begin{equation}
\label{Ball eN+1}
\langle e_0, e_{N+1}\rangle=-\frac{i}{2},\ \langle e_\alpha, e_{N+1}\rangle=0,\
\langle e_{N+1}, e_{N+1}\rangle=0.
\end{equation}

\bigskip

\section{Relationship among four definitions of $II_M$}

\begin{lemma}
\label{3 II M}
Let $H: M' \ra \p \HH^{N+1}$ be a CR smooth embedding
where $M'$ is a strictly pseudoconvex smooth real hypersurface in $\CC^{n+1}$.
We denote $M=H(M')$. Then the following statements are equivalent:

(i) The CR second fundamental form $II_M$ by Definition 1 identically vanishes.

(ii) The CR second fundamental form $II_M$ by Definition 2 identically vanishes.

(iii) The CR second fundamental form $II_M$ by Definition 3 identically vanishes.

(iv) The CR second fundamental form $II_M$ by Definition 4 identically vanishes.
\end{lemma}

\bigskip

\noindent{\it Proof}\ \ \ \ (i) $\Longleftrightarrow$ (ii)  by (\ref{IIM b}).

(iii) $\Longleftrightarrow$ (iv)  \ The equivalence follows by the facts that, for Definition 3 and 4,
$II_M^e\equiv 0$ for one first-order adapted lift $e$ if
and only if $II^s_M \equiv 0$ for any first-order adapted lift $s$,  that a first-order adapted lift from $M$ to $SU(N+1,1)$ must be a first-order adapted lift from $M$ to $GL^Q(\CC^{N+2})$.

(iv) $\Longrightarrow$ (i): Let $M\subset \p\HH^{N+1}$ be a $(2n+1)$ dimensional CR submanifold with
CR dimension $n$ that admits a first-order adapted lift $e$ into $SU(N+1,1)$.
Consider the pull-backed Maurer-Cartan form over $M$ by $e$
\[
\omega
=\begin{pmatrix}
\omega^0_0 & \omega^0_\beta & \omega^0_\nu & \omega^0_{N+1}\\
\omega^\alpha_0 & \omega^\alpha_\beta & \omega^\alpha_\nu &
\omega^\alpha_{N+1}
\\
0 & \omega^\mu_\beta & \omega^\mu_\nu& \omega^\mu_{N+1}\\
\omega^{N+1}_0 & \omega^{N+1}_\beta & 0&
\omega^{N+1}_{N+1}\\
\end{pmatrix}.
\]
with
\begin{equation}
\label{3 basic formula theta 3}
\begin{array}{lll}
&& \omega^0_0 + \ov{\omega^{N+1}_{N+1}} =0,\
\omega^{N+1}_0=\ov{{\omega}_0^{N+1}},\
\omega^0_{N+1}=\ov{\omega_{N+1}^0},\\
&&
\omega^{N+1}_A=  2i \ov{\omega^A_0},\ \omega^A_{N+1}=- \frac{i}{2}
\ov{\omega^0_A},\
\omega^A_B+\ov{\omega^B_A}=0,\ \omega^0_0+\omega^A_A+\omega^{N+1}_{N+1}= 0,
\end{array}
\end{equation}
where $1\le A\le N$.

Let $\theta=\omega^{N+1}_0$ which is a real $1$-form by (\ref{3 basic formula theta 3}).
By $d\omega=-\omega\wedge \omega$ and (\ref{3 basic formula theta 3}), we obtain
\[
d\theta=-\omega^{N+1}_0\wedge \omega^0_0-\omega^{N+1}_\alpha\wedge \omega^\alpha_0
-\omega^{N+1}_{N+1}\wedge \omega^{N+1}_0
=2i \omega^\alpha_0\wedge\ov{\omega^\alpha_0}-\theta\wedge (\omega^0_0+\ov{\omega^0_0})
=i\theta^\alpha\wedge \ov{\theta^\alpha},
\]
where we denote
\begin{equation}
\label{Definition 4 to 1}
\theta^\alpha=\sqrt 2 \omega^\alpha_0+c_\alpha \theta
\end{equation}
for some functions $c_\alpha$. Therefore, (\ref{orthonormal}) holds and hence $M$ is a strictly
pseudoconvex pseudohermitian manifold
with an admissible coframe $(\theta, \theta^\alpha)$. Hence Definition 4
of $II_M\equiv 0$ implies Definition 1 of $II_M\equiv 0$.

\bigskip

(i) $\Longrightarrow$ (iv): Definition 1 of $II_M$ gives a
coframe $(\theta, \theta^\alpha)$ which corresponds to Definition 2 of $II_M$ with respect
to a defining function $\rho$ of $M$ in $\p\HH^{N+1}$.

Now take a first-order adapted lift $e$ from $M$ into $SU(N+1,1)$.
By (\ref{Definition 4 to 1}), it corresponds to a coframe $(\theta, \theta^\alpha)$ on $M$ and
by (\ref{bilinear form 2}), it corresponds Definition 2 of $II_M$ by some choice of the defining
function $\hat \rho$ of $M$ in $\p\HH^{N+1}$.

The above $\rho$ and $\hat \rho$ may not be the same. But Definition 2 of $II_M\equiv 0$ is independent of choice
of defining functions, which gives (i) $\Longrightarrow$ (iv). \ \ \ $\Box$

\section{Proof of Theorem 1.1}
\begin{lemma}
\label{spherical lemma}
(cf. [EHZ04], corollary 5.5)
Let $H: M' \ra M \hookrightarrow \p\HH^{N+1}$ be a smooth CR embedding of a strictly pseudoconvex smooth real hypersurface $M \subset \CC^{n+1}$. Denote by $(\omega^{\ \mu}_{\alpha\ \beta})$ the CR second
fundamental form matrix of $H$ relative to an admissible coframe
$(\theta, \theta^A)$ on $\p\HH^{N+1}$ adapted to $M$.
If $\omega^{\ \mu}_{\alpha\ \beta} \equiv 0$ for all $\alpha, \beta$ and $\mu$, then $M'$ is locally
CR-equivalent to $\p\HH^{n+1}$.
\end{lemma}

\bigskip

\noindent{\it Proof of Theorem 1.1}\ \ \ \ {\bf Step 1. Reduction to a problem for
geometric rank}\ \ \ \ By Lemma \ref{3 II M} and Lemma
\ref{spherical lemma} and the
hypothesis that the CR second fundamental form identically vanishes,
we know that $M$ is locally CR equivalent to $\p\HH^{n+1}$.

Then $M$ is the image of a local smooth CR map $F: U \subset \p\HH^{n+1} \ra M
\subset \p\HH^{N+1}$ where $U$ is a open set in $\p\HH^{n+1}$ .
By a result of Forstneric\cite{Fo89}, the map $F$ must be a rational map.
It suffices to prove that
$F$ is equivalent to a linear map.
By Lemma \ref{rank 0},
it is sufficient to prove that the geometric rank of $F$ is zero: $\kappa_0=0$.

Suppose $\kappa_0>0$ and we seek a contradiction.

\bigskip

\noindent{\bf Step 2. Reduction to a lift of  $\big((H\circ \tau^F_p)(M), 0\big)$}\ \ \ \
Take any point $p\in U \subset \p\HH^{n+1}$ with $\kappa_0=\kappa_0(p)>0$, and consider the associated map
(see Lemma \ref{**normalform})
\begin{equation}
\label{F***}
F_p^{***}=H\circ \tau^F_p \circ F\circ \sigma^0_p\circ G: \p\HH^{n+1} \ra
\p\HH^{N+1}, \ \ \ F^{***}_p(0)=0,
\end{equation}
where $\sigma^0_p$ is defined in (\ref{sigma}), $\tau^F_p$ is defined in (\ref{tau}),
$G\in Aut_0(\HH^{n+1})$ and $H\in Aut_0(\HH^{N+1})$ are automorphisms.
By Theorem \ref{***normalform},
$F^{***}_p=(f,\phi,g)$ satisfies the following normalization conditions:

\begin{equation}
\label{proof of them 1}
\left\{
\begin{aligned}
f_j=&z_j+\frac{i\mu_j}{2}z_jw+o_{wt}(3),\ \ \frac{\p^2 f_j}{\p w^2}(0)=0,\
j=1
\cdots,\kappa_0,\ \mu_j>0,\cr
f_j=&z_j+o_{wt}(3), \ \ j=\kappa_0+1,\cdots, n-1\cr
             g = & w+o_{wt}(4),\cr
\phi_{jl}=&\mu_{jl}z_jz_l+o_{wt}(2), \ \h{where}\ \ (j,l)\in {\cal S}\
\h{with}\ \ \mu_{jl}>0 \ \h{for}\ \  (j,l)\in {\cal S}_0\cr
 &\h{and}\ \mu_{jl}=0\ \h{otherwise}
\end{aligned}
\right.
\end{equation}
where $\mu_{jl}=\sqrt{\mu_j+\mu_l}$ for $j,l\le
\kappa_0$ $j\not=l$,
$\mu_{jl}=\sqrt {\mu_j}$ if $j\le \kappa_0$ and $l>\kappa_0$ or if
$j=l\le \kappa_0$. Here the assumption that $\kappa_0>0$ is used.

From (\ref{F***}) we obtain
\[
\begin{matrix}
(M, F(p)) & \xrightarrow{H\circ \tau^F_p} & \big(H\circ \tau^F_p(M), 0\big)\\
\uparrow F &\ & \uparrow F^{***}_p \\
(\p\HH^{n+1}, p) & \xleftarrow{\sigma^0_p\circ G} & (\p \HH^{n+1}, 0) \\
\end{matrix}
\]
If we can show that there exists a first-order adapted lift $e$ from the submanifold
$H\circ \tau^F_p(M)$ near $0$ into $SU(N+1,1)$ such that the corresponding
CR second fundamental form
\begin{equation}
\label{reduction 1}
II^e_{H\circ \tau^F_p(M)}\not= 0 \ at\ 0,
\end{equation}
then we obtain a first-order adapted lift $\w e:=\big(H\circ \tau^F_p \big)^{-1} \circ e \circ
H\circ \tau^F_p$ from the submanifold
$M$ near $F(p)$ into $GL^Q(\CC^{N+1})$ such that the corresponding
CR second fundamental form
\begin{equation}
\label{reduction 11}
II^{\w e}_{M}\not= 0 \ at\ F(p).
\end{equation}
Notice that the map $H\circ \tau^F_p \in GL^Q(\CC^{N+2})$ but $H\circ \tau^F_p \notin
SU(N+1, 1)$, so that the lift $\w e$ is not from $M$ into $SU(N+1,1)$. This is why we
have to introduce Definition 3.

Since we take arbitrary $p\in \p\HH^{n+1}$, from (\ref{reduction 11}) it
concludes that $II_M\not\equiv 0$, but this is a desired contradiction.

\bigskip

\noindent{\bf Step 3. Calculation of the second fundamental form}\ \ \ \
It remains to prove existence of the lift $e$ such that (\ref{reduction 1}) holds.

The lift $e$ constructed in the second half of Section 7 is a first-order adapted lift from
$H\circ \tau^F_p(M)$ near $0$ into $SU(N+1,1)$ which defines a CR second
fundamental form as a tensor $II^e_{H\circ \tau^F_p(M)}=q^\mu_{\alpha \beta}
\omega^\alpha \omega^\beta\otimes (\underline{e_\mu})$ in
(\ref{II M 4 2}). If we can show
\begin{equation}
\label{reduction 2}
q^\mu_{\alpha\beta}(0)=\frac{\p^2 \phi_\mu}{\p z_\alpha \p z_\beta}\bigg|_0,
\end{equation}
where $F^{***}_p=(f, \phi, g)=(f_\alpha, \phi_\mu, g)$. Since we assume that
$\kappa_0>0$, by (\ref{proof of them 1}) and (\ref{reduction 2}), it implies
$q^\mu_{\alpha\beta}(0)\not=0, \forall \alpha, \beta$ and $\mu$, i.e.,
$II^e_{H\circ \tau^F_p(M)}\not=0$. This proves (\ref{reduction 1}).


Let $E=(e_0, e_\alpha, \hat E_\mu, e_{N+1})$ be the lift constructed in Theorem \ref{existence of lift}
(see the remark at the end of the proof of Theorem \ref{existence of lift}) and
in (\ref{Ball e0}) (\ref{Ball ealpha}) and (\ref{Ball eN+1}).
Since $E|_0=Id$, we have
\[
\omega|_0=(E^{-1}|_0) (d E)|_0=d E|_0
\]
so that
\[
\omega|_0=\begin{bmatrix}
0 & * & ... & *\\
dz_1 & * & ... & *\\
\vdots & \vdots& \  &\vdots  \\
d z_n & * & ... & *\\
* & * & ... & * \\
\vdots & \vdots & \  &\vdots  \\
* & * & ... & * \\
dw & * & ... & *\\
\end{bmatrix} \bigg|_0.
\]
Hence $\omega^1_0|_0=dz_1,\ ...,\ \omega^n_0|_0=d z_n,\ \omega^{N+1}_0|_0=d z_{N+1}$.
Then by applying the chain rule, we obtain
\[
\omega^\mu_j|_0=d E^\mu_j|_0 = d \big(( L_j \phi_\mu)\circ (F_{fg})^{-1} \big)\big|_0
= \frac{\p}{\p z_k}\big( ( L_j \phi_\mu)\circ (F_{fg})^{-1}  \big)|_0
d z_k=\frac{\p^2 \phi_\mu}{\p z_k \p z_j}|_0 \omega^k_0|_0,
\]
for any $j,k\in\{1,2,...,n,N+1\}, \ n+1\le \mu\le N$. Hence (\ref{reduction 2}) is proved.
The proof of Theorem 1.1 is complete. \ \ \ $\Box$

\bigskip
\bibliographystyle{amsalpha}

Shanyu Ji (shanyuji$@$math.uh.edu), Department of Mathematics, University of Houston, Houston, TX 77204;

Yuan Yuan (yuanyuan$@$math.rutgers.edu), Department of
Mathematics, Rutgers University, Piscataway, NJ 08854.

\end{document}